\documentclass[12pt]{article}

\usepackage{amsfonts,amssymb,amsmath}

\usepackage{epsfig,graphics,graphicx,color}
\usepackage{psfrag}
\newcommand{\rd}{\rm d}
\newcommand{\wt}{\tilde}
\newcommand{\wh}[1]{\widehat{#1}}

\newcommand{\ovl}[1]{\overline{#1}}
\newcommand{\reff}[1]{(\ref{#1})}

\newcommand{\E}{{{\bf {E}}}}

\newcommand{\proofr}{{\bf Proof. }}

\newtheorem{theorem}{Theorem }
\newtheorem{lemma}{Lemma }
\newtheorem{utv}{Preposition }

\begin{document}

\author{E.A. Pechersky and N.D. Vvedenskaya \\
\footnotesize{ Dobrushin Mathematical Laboratory,}\\
\footnotesize{Institute for Problems of
Information Transmission, Russian Academy of Sciences,}\\
\footnotesize{B.~Karetnyi per., GSP--4 Moscow 101447, Russia}\\
\footnotesize{pech@iitp.ru, ndv@iitp.ru}}
\title{A circle of interacting servers; spontaneous collective
behavior in case of large fluctuations}

\date{}

\maketitle

\begin{abstract}
We consider large fluctuations, namely   overload of servers, in a
network with dynamic routing of messages. The servers form a
circle. The number of input flows is  equal to the number  of
servers, the messages of any flow are distributed between two
neighboring servers, upon its arrival a message is directed to
the least loaded of these servers. Under the condition that at
least two servers are overloaded the number of overloaded servers
in such  network depends on the rate of input flows. In
particular there exists critical level of input rate that in case
of higher rate most probable that all servers are overloaded.
\end{abstract}

\bigskip

\section{Introduction}
\medskip

This work presents an effect in a network with interacting servers
that  can be called a spontaneous collective behavior in case of
large fluctuations.

We consider networks with dynamic routing of messages. In such
networks the server to which a message is directed depends on
the  network's state at the message arrival moment. One of the
problems arising here is the analysis of probability of large
fluctuations, for example probability of large delays.

There are many works where  large fluctuations in  networks with
dynamic routing have been investigated. In \cite{AH} - \cite{DPSV}
the networks with two servers and three independent input flows
have been considered where only one flow is divided between two
servers depending either on the workload of servers or on the
queue lengths. In  \cite{PV}  a network with a group of  servers
and several flows has been considered where each flow is assigned
to some subgroup of servers; upon its arrival a message selects a
server with the shortest queue(i.e. a queue with least number of
messages). In this work the large deviation principal for the
flows upon the servers is proved. We would like to stress that
the service time of messages in \cite{PV} is exponentially
distributed. That allows to use Markov property for the flows to
servers (after splitting of the input flows.)

The more full list of references can be found in the mentioned
works.


\begin{figure}
{\includegraphics[natheight=9cm,natwidth=3cm,scale=1.0]{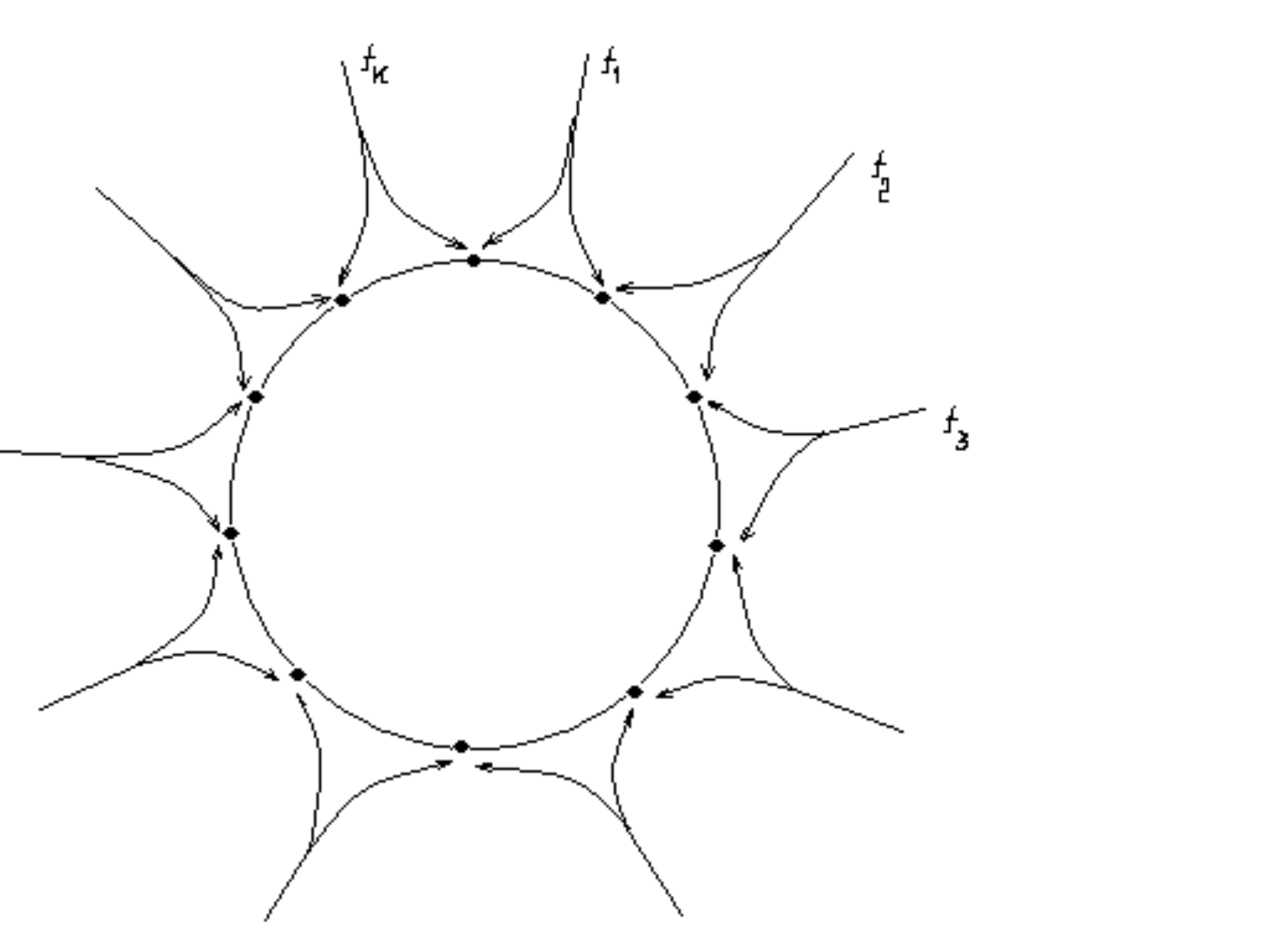}}
\caption{}
\end{figure}

Here we consider circle networks that are formed by $k$ servers
and $k$ identical independent  input Poisson   flows. Messages of
any flow are assigned to two nearby servers (see Figure 1) A
message direction  depends on the  workloads on these two
servers, namely upon its arrival a message is directed to the
server with the smallest workload. The constant  speed of work of
each server is equal to 1, the discipline is "first in - first
out" (FIFO). If a message finds the server busy it is put into a
infinite buffer to wait for service. We consider the networks
that work stationary. That means that  with probability 1 the
queues do not increase infinitely. But there may appear  large
fluctuations, for example during a short period one of  flows may
bring  very large amount of work. We  say that during this period
the flow is overheated. Suppose for example that the flow $f_{1}$
is overheated. Then the buffers of servers $s_1$ and $s_k$ that
are assigned to this flow will contain a large amount of work.
What is  the behavior of other flows ?

 In this paper we  show that in case where $f_1$ is overheated there
exist at least two scenarios of network performance. What scenario
is more probable and  is realized depends on the rate of input
flows. Namely, we show that there exists such  value of input flow
rate that in case when the arrival rates are above it the
overheating of flow $f_1$ coincides with overheating of all flows.
Such behavior may be regarded as the spontaneous appearance of
collective behaviour. On the other hand  in case of low rate the
overheating of $f_1$ does not bring the overheating of other
flows. There may exist also intermediate stages where the
overheating of one flow coincides with the overheating of a
number of neighboring flows. The existence of such intermediate
stages depends on  message length distribution of input flows.
For some length distributions,  for example for exponential one,
the intermediate stages do not exist.

The proof of this result is based on large deviation principle of
\cite{DP} (see also \cite{LP}). We do not present the details of
this reduction.  The details of  application of large deviation
principle for several other problems  one can fined  in \cite{DP,
PSV1, AP}. Instead of the proof we present some ideas needed for
the   proof and  some fragments of rigorous arguments. Therefore
the following text is not   a proof in the conventional sense.

In the next section we introduce several  concepts and definitions
and describe the main ideas of our approach.  In section 3 several
examples are presented.

\medskip

\section{Cyclic networks. Main result}

\medskip

The cyclic network of size $k,\ k\geq 3,$ is formed by $k$ servers
and $k$ input flows. Let $F=\{f_1,...,f_k\}$ be the set of flows and
$S=\{s_1,...,s_k\}$  the set of servers. All servers work with speed
1, the  discipline is FIFO. Each server has an infinite buffer. The
messages of flow $f_i\in F$ are served by two servers $s_i\in S$
 and $s_{i-1}\in S$ (here  always $i\mp 1\mod k$). Each message of $f_i$ is
 directed to that of  servers $s_i,\ s_{i-1}$ which  becomes  idle first.

Let the random sequence $(\xi_{n}^{i},\tau_{n}^{i})$,
$n=...,-1,0,1,...$, describes the  input flow $f_i\in F$. The
random variables $\tau_{n}^{i}$  are the intervals between  the
arrivals of  messages of flow $f_i$. The random variables
$\xi_{n}^{i}$ are the lengths of messages. We consider
homogeneous Poisson flows therefore all
$(\xi_{n}^{i},\tau_{n}^{i})$ are independent and equally
distributed,  $\tau_{n}^{i}$ are exponentially distributed with
rate $\lambda$, i.e. $\Pr(\tau_{n}^{i}>t)=e^{-\lambda t}$. The
distributions of input flow do not depend on $i$.

Introduce a random pair $(\xi,\tau)$ that is distributed as any pair
$(\xi_{n}^{i},\tau_{n}^{i})$. Each sequence
$(\xi_{n}^{i},\tau_{n}^{i})$ is numbered so that $\tau_{0}^{i}<0$ and
 $\tau_{1}^{i}\geq 0$.

We propose that  $\xi$ has exponential moments, i.e. there exist
such $\theta_{+}$, \;\;$0<\theta_{+}\leq \infty$ that
\begin{equation*}
\label{lapl} \varphi(\theta)=\E e^{\theta\xi}<\infty\ \ \mbox{
for}
 \ \  \theta<\theta_{+},
\end{equation*}
\begin{equation}\label{phi+}
\lim_{\theta\uparrow\theta_{+}}\varphi(\theta)=\infty.
\end{equation}
The condition that guarantees the  stationarity and even ergodicity
is
\begin{equation}\label{ergodic}
\lambda\varphi^{\prime}(0)<1.
\end{equation}
     This condition is intuitively obvious but its sufficiency
 for the existence of  stationary stage needs a proof. In
particular one may construct a Liapunov function that shows that
with large probability the workloads of all stations  are located in
a compact region. We omit this construction.

\medskip

To formulate the problem that is investigated we introduce a
notion of {\it virtual message}. That is a message that arrives
with flow $f_i$ at a nonrandom time moment, for example at 0 and
has zero length. It  joins a queue at $s_i$ or $s_{i-1}$,
depending on where the workload is less, and waits there its
service. Let $\omega_{i}$ be the $i$th virtual message's waiting
time. We are interested in the probability of large waiting for
virtual massage that arrived with $f_1$, i.e. $\Pr (\omega_1\geq
d)$, where $d$ is large. In some cases, for special choice of
$\xi$ distribution the needed probability can be calculated. But
in general case it is impossible to present the explicit
expression. Therefore one looks for the asymptotics of probability
\begin{equation}\label{virt}
\Pr(\omega_{1}>nd),
\end{equation}
where $ n\to\infty$. More detailed , one looks for
\begin{equation}\label{I}
J=\lim_{n\to \infty}\frac{-1}{n}\ln \Pr(\omega_{1}\geq nd).
\end{equation}
This problem belongs to the theory of large deviation. We  use this
theory even where it is possible to get an explicit solution.

\medskip

\subsection{Collective fluctuations}

\medskip

Next theorem states the existence of input flow values
$\lambda_k\leq \lambda^k$ such that separate different types of
system performance. Let
$$ \hat{\lambda}=\frac{1}{\varphi'(0)}.
$$
As we mentioned before the system performance is stationary if
$\lambda<\hat{\lambda}$.

We say that  the flow $f_i$ is overheated at $t=0$ if
$\Gamma_{i}(n)=(\omega_{i}\geq nd)$. It means that  a virtual
message that arrived with $f_i$ is waiting for being served  at
least $nd$.
 Let $\ovl{\Gamma}_{i}(n)$ be the complement to event
 $\Gamma_{i}(n)$.

\medskip

\begin{theorem}.\label{t}
For any system size $k,\;k\geq 3$

and for any $d$ there exist   $\lambda_{k},\,\lambda^{k}$ (not
depending on $d$), such that
$0<\lambda_{k}\leq\lambda^{k}<\wh{\lambda}$,
 and
\begin{enumerate}
\item [${\blacktriangleright}$]\ if $\lambda<\lambda_{k}$
then
\begin{enumerate}
\item [1)]
$\Pr\left(\Gamma_{1}(n)\bigcap_{i=2}^{k}\ovl{\Gamma}_{i}(n)\,
\big| \Gamma_{1}(n)\right)\to 1$, as $n\to \infty$.\\

\item [2)] $J=\lim_{n\to\infty}\frac{-1}{n}\ln
\Pr\left(\Gamma_{1}(n)\bigcap_{i=2}^{k}\ovl{\Gamma}_{i}(n)\right)
=2\theta(\lambda ,1)d$\\
where $\theta(\lambda,1)$ is a positive root to equation
\begin{equation}\label{theta1}
 2\theta=\lambda[\varphi(\theta)-1],
\end{equation}
\end{enumerate}
\item [${\blacktriangleright}$]\ if $\lambda>\lambda^{k}$
then
\begin{enumerate}
  \item [3)] $\Pr\left(
\bigcap_{i=1}^{k}\Gamma_{i}(n)\, \big| \Gamma_{1}(n)\right) \to 1$
as $n\to \infty$,
\item [4)]$J=\lim_{n\to\infty}\frac{-1}{n}\ln
\Pr\left(\bigcap_{i=1}^{k}\Gamma_{i}(n) \right) =k\theta^{*}d$,\\
where $\theta^{*}$ is a positive root to equation
\begin{equation}\label{thetak}
\theta=\lambda[\varphi(\theta)-1].
\end{equation}
\end{enumerate}
\end{enumerate}
\end{theorem}
It follows form this statement  that in case $\lambda<\lambda_{k}$
the main contribution to the  probability of event
$\Gamma_1(n)=(\omega_{1}>nd)$ brings the flow $f_1$  which  the
virtual message arrives with. Only this flow is overheated, the
others stay not overheated. But in case $\lambda>\lambda^{k}$ all
input flows are overheated. Though  the virtual message is
combined with flow $f_1$ only its large delay is connected with
an effect of collective behavior that is similar to the effect of
spontaneous magnetization
 in  statistical physics systems.

It is possible to present a more   accurate statement. In the next
section we introduce several definitions and present a precise
theorem.
\medskip

\subsection{Random processes connected with the system and large deviations }

\medskip

Below we introduce several random processes connected with the
cyclic network.

\begin{enumerate}

  \item \label{zeta}Process

$$\zeta^{i}(t)=
        \sum\limits_{j:\ \sum_{r=1}^{j}\tau_r^i<t}\xi_j^i,
$$
describes the amount of work brought by $f_i$ during time interval
$[0,t)$.

  \item \label{wi} $k$ dimensional Markov process $w^{1},...,w^{k} $
defined by  generator
\begin{align}
&Lg(w^{(1)},...,w^{(k)})=\nonumber\\ &\hphantom{aaa}=\
-\sum_{i=1}^{k}\frac{\partial g}{\partial
w^{(i)}}{\mathbb{I}}_{(w^{(i)}>0)}(w^{(i)})+\nonumber\\
&\hphantom{aaa}+\lambda\sum_{i=1}^{k}\left[ {\E}
g(w^{1},...,w^{i}+\xi,...,w^{k})-\right.\nonumber\\
&\hphantom{aaaaaaaaaaaa}\left.-g(w^{1},...,w^{i},..., w^{k})\right]
{\mathbb{I}}_{(w^{i}\leq w^{i+1})}(w^{i},w^{i+1})+ \nonumber\\
&\hphantom{aaa}+\lambda\sum_{i=1}^{k}\left[ {\E}
g(w^{(1)},...,w^{i+1}+\xi,...,w^{k})-\right.\nonumber\\
&\hphantom{aaaaaaaaaaaa}\left.-g(w^{1},...,w^{i+1},...,
w^{k})\right]{\mathbb{I}}_{(w^{i}>w^{i+1})}(w^{i},w^{i+1}),
\nonumber
\end{align}
 describes the loading of all servers. We propose that
$w^{i}(0)=0,\;i=1,...,k$. The symbol $\mathbb{I}_{A}$ is the
indicator of  a set $A$,  $g$  is arbitrary  differentiable
function in $L_2(\mathbb{R}^{k})$. The second and third sums of
the right hand side indicate that a message of length $\xi$
arriving with $f_{i+1}$ is directed to one of
 servers $s_i$, $s_{i+1}$ depending on where the loads $w^i$, $w^{i+1}$
 is less.

 A typical trajectory of this process is a $k$ dimensional function
  with non negative components. Each component is  formed by jumps and
 piecewise linear functions between the jumps. On the intervals
 where the trajectory is differentiable its derivative is either
 $-1$  (if it is positive) or 0 (if it  equals 0).

\item\label{whwi}
Process $(\wh{w}^{1}(t),...,\wh{w}^{k}(t))$  describes the amount
of work that was brought to each of $k$  servers presuming that
${w}^{i}(0)=0$. Let $t_r^i>0$, $r=1,...$ be the set of jumps of
process $w^i$ and $\sigma^i_r$, $r=1,...$  the sizes of these
jumps. Then
\begin{equation}\label{hw}
\wh{w}^{i}(t)=\sum_{r:\:t_{r}^{i}<t}\sigma_{r}^{i}.
\end{equation}

The trajectories of this process are  step functions, where jumps
coincide in size and time  with those of  process $w^{1},...,w^{k}
$.
\end{enumerate}

We define these three $k$-dimensional processes
$\zeta=(\zeta^{1},...,\zeta^{k})$, $w=(w^{1},...,w^{k})$ and
$\wh{w}=(\wh{w}^{1},...,\wh{w}^{k})$ onto the same probability
space. Here the jumps of $w$ and $\wh{w}$ coincide in size and time.
Let $\zeta^i$ have a jump at time moment $t$. Then at the same
moment either $w^{i-1}$ or $w^{i}$ has a jump of the same size,
depending on whether an inequality $w^{i-1}(t-0)<w^{i}(t-0)$ or the
opposite one takes place. In case of equality we propose that
process $w^{i-1}$ has a jump.

Defining the processes on one probability space one can define a
mapping
\begin{equation}\label{jumpo}
G:\:\mathcal{X}^{o}\to\mathcal{X}^{o}
\end{equation}
of the set realization of processes $\zeta$ onto the set realization
of processes $\wh{w}$.  Here $\mathcal{X}^{o}$ is the set of
nondecreasing stepwise functions that are equal to $0$ at $0$.

All event that are considered here are connected with asymptotic
characteristics of  network as $n\to\infty$. Therefore we
introduce the scaled versions of processes and events.

We start  with introduction of scaled versions of the processes
\ref{zeta} and  \ref{whwi}:
\begin{enumerate}\label{processn}
  \item \label{zetan}
$$\zeta^{i}_{n}(t)= \frac{1}{n}\zeta^{i}(nt).
$$
\item\label{whwin}
\begin{equation*}
\wh{w}^{i}_{n}(t)=\frac{1}{n}\wh{w}^{i}(nt).
\end{equation*}
\end{enumerate}

Being interested in  asymptotic behavior of probabilities of
different events we introduce also several notions defined by the
sequences of events.

Let $\{x_{T,a}(\cdot)\}$ be a piecewise  linear trajectories
\begin{equation}\label{xx}
\dot x_{T,a}(t)=\begin{cases}
a&\mbox{ if }t\in [0,T]\\
\lambda\varphi^{\prime}(0)&\mbox{ if }t>T
\end{cases}
\end{equation}

We say that a flow $f_i$ is overheated on an interval $[0,T]$ if
following event takes place
\begin{equation*}
\mathcal{A}_{\varepsilon,n}^{i}(x)=
\left(\sup_{t\in[0,T]}\left|x(t)-\zeta_{n}^{i}(t)
\right|<\varepsilon \right),\;\;x=x_{T,a},\;\;a>1, \varepsilon>0.
\end{equation*}
Remark that the notion of overheated flow is introduced only for
piecewise linear trajectories. We do not need more general
definition.

We say that a flow $f_i$ is not overheated on an interval $(0,T)$ in
case of event
 $\mathcal{A}_{\varepsilon,n}^{i}(x)=\left(\sup_{t\in[0,T]}\left|x(t)-\zeta_{n}^{i}(t)
\right|<\varepsilon\right),\;\;x=x_{T,a},\;\;a\leq
1,\;\;\varepsilon>0. $ If $a=a^{(0)}=\lambda\varphi^{'}(0) $ then
not overheated input flow is close to its mean value.

\medskip

Remark that the distribution of  $\omega_{1}$ coincides with the
distribution of $\min\left\{\sup_{t\geq 0
}\{\wh{w}^{1}(t)-t\},\;\sup_{t\geq 0 }\{\wh{w}^{k}(t)-t\}\right\}$.
That can be shown easily, we omit the needed construction. (For the
problems considered in \cite{DP, PSV1, AP} similar facts are
explained in these works, the needed construction for a one-channel
systems one can also find in monograph \cite{B} ).

The probability \reff{virt} can be expressed in terms of scaled
processes $\wh{w}^{i}_{n}(t)$. It is equal to
\begin{equation}\label{virt1}
\Pr(\omega_{1}>nd)=\Pr\left(\wh{\mathcal{M}}_{n}\right),
\end{equation}
where event
\begin{equation*}
\wh{\mathcal{M}}_{n}=\left(\min\left\{\sup_{t\geq 0
}\Big\{\wh{w}^{1}_{n}(t)-t\Big\},\;\sup_{t\geq 0
}\Big\{\wh{w}^{k}_{n}(t)-t\Big\}\right\}>d\right)
\end{equation*}
means that both processes $\wh{w}^{1}_{n}$ ш $\wh{w}^{k}_{n}$
intersect the line $u(t)=d+t$.

 Using the mapping $G$ of process $\zeta_{n}$ realizations onto  process
 $\wh{w}_{n}$ realizations   (see \reff{jumpo}) we can introdus an
 event $\mathcal{M}_{n}=G^{-1}\wh{\mathcal{M}}_{n}$.

\medskip

\begin{theorem}.\label{t1}
For any size $k,\;k\geq 3$ of the network and for any $d>0$ there
exist not depending on $d$ \ \ $\lambda_{k}$ and $\lambda^{k}$,
$0<\lambda_{k}\leq\lambda^{k}<\wh{\lambda}$, such that
\begin{enumerate}
\item [${\blacktriangleright}$]\ if $\lambda<\lambda_{k}$ then there exist
 $T^{*}_{1}$ and $a^{*}_{1}$ such that
\begin{enumerate}
\item [1)]$ \lim_{\varepsilon \to 0}\lim_{n\to 0}
\Pr\left(\mathcal{A}_{\varepsilon,n}^{1}(x_{T^{*}_{1},a^{*}_{1}})
\bigcap_{i=2}^{k}{\mathcal{A}}_{\varepsilon,n}^{i}({x}_{T^{*}_{1},a^{(0)}})\Big|\mathcal{M}_{n}\right)
\to 1$
\item [2)] $J=\lim_{\varepsilon \to 0}\lim_{n\to 0}
\frac{-1}{n}\ln\Pr\left(\mathcal{A}_{\varepsilon,n}^{1}(x_{T^{*}_{1},a^{*}_{1}})
\bigcap_{i=2}^{k}{\mathcal{A}}_{\varepsilon,n}^{i}({x}_{T^{*}_{1},a^{(0)}})\right)
=2\theta(\lambda ,1)d$,\\
where $\theta(\lambda,1)$ is the positive root to equation
\begin{equation}\label{theta2}
 2\theta=\lambda[\varphi(\theta)-1],
\end{equation}
\item[3)] besides
\begin{equation}\label{at1}
a^{*}_{1}=\lambda\varphi^{'}(\theta(\lambda ,1)),\;\;\;
   T^{*}_{1}= \frac{d}{\lambda\varphi^{'}(\theta(\lambda,1))-1}.
\end{equation}
\end{enumerate}
\item [${\blacktriangleright}$]\ if $\lambda>\lambda^{k}$
then there exist  $T^{*}_{k}$ and $a^{*}_{k}$ such that
\begin{enumerate}
  \item [4)] $\lim_{\varepsilon \to 0}\lim_{n\to 0}\Pr\left(
\bigcap_{i=1}^{k}{\mathcal{A}}_{\varepsilon,n}^{i}(x_{T^{*}_{k},a^{*}_{k}})\Big|\mathcal{M}_{n}\right)
\to 1$
\item [5)]$J=\lim_{\varepsilon \to 0}\lim_{n\to 0}\frac{-1}{n}\ln
\Pr\left(
\bigcap_{i=1}^{k}{\mathcal{A}_{\varepsilon,n}^{i}(x_{T^{*}_{k},a^{*}_{k}}})
 \right) =k\theta^{*}d$,\\ where $\theta^{*}$ is the positive root to equation
\begin{equation}\label{thetak1}
\theta=\lambda[\varphi(\theta)-1],
\end{equation}
\item[6)] besides
\begin{equation}\label{at2}
a^{*}_{k}=\lambda\varphi^{'}(\theta^{*}),\;\;
   T^{*}_{k}= \frac{d}{\lambda\varphi^{'}(\theta^{*})-1}.
\end{equation}

\end{enumerate}
\end{enumerate}
\end{theorem}
\medskip

This more precise form of  the theorem presents the mean dynamic
of conditional process under the condition of event
$\mathcal{M}_{n}$. Namely, the trajectory
$(x_{T^{*}_{1},a^*_{1}},a^{(0)}t,...,a^{(0)}t)$ is the conditional
mean dynamic of processes $(\zeta^{1}_{n},...,\zeta^{k}_{n})$ on
the interval $[0,T^{*}_{1}]$ in case $\lambda<\lambda_{k}$, and
the trajectory
$(x_{T^{*}_{k},a^*_{k}},...,x_{T^{*}_{k},a^{*}_{k}}) $ is the
conditional mean dynamic
 on the interval $[0,T^{*}_{k}]$ in case $\lambda>\lambda^{k}$

\medskip

\subsection {Ideas of proof}
\medskip

 The value of   $J$ could  be found if the
large deviation principle for processes $\wh{w}^{1}(t)$ and
$\wh{w}^{k}(t)$ would be known. These processes are not Poisson.
But they are the functionals of   ${\zeta}^{1}(t),...\zeta^k(t)$.
We remind that they are defined on the same probability space.
The processes $\zeta^i, i=1,...,k$  are  Poisson. Therefore one
can use some known results on large deviations of Poisson
processes. We use the results on large deviation principle from
\cite{DP} (see also \cite{LP}).

For the problems that are considered here the application of large
deviation principle consists of two parts. First one has to check
the validity of the large deviation principle and describe the
corresponding to the problem event. After that one has to find
minimum of rate function on the event. As usual,  application of
large deviation principle reduces  the problem   to the search of
the point (sometimes several points) of the event where the rate
function take minimal values. A small neighborhood of the point
of the minimum brings the main contribution to the asymptotic of
logarithm of event probability. In our case the event is some set
of trajectories and the rate function is an integral functional
on the trajectories. Therefore a search of trajectory that
minimizes the rate function is reduced to  a variational problem.

To check the validity of large deviation principle one has to
present the topological space, the sequence of measures on this
space that tend to a $\delta$-measure and the rate function.

The topological space is a set of non decreasing $k$ dimensional
functions on $[0,\infty)$
\[\wh{\mathcal{X}}=\mathfrak{X}^{k}=
\{(x^{1}(\cdot),...,x^{k}(\cdot))\},\] equipped with {\it
uniformly  week} topology.

 The detailed description of these functions one can find in
 \cite{DP}, where the uniformly  week topology is introduced on this
 space. We use this topology. It is explained in the same paper why the
 week topology is too week and one has to use the uniformly  week
 topology for the problems of the considered kind.

The sequence of measures $P_{n}$ is defined by processes
$\zeta_{n}$.

Finally, the rate function $\wh{I}$ is defined on
$\langle\wh{\mathcal{X}},P_{n}\rangle$ by equalities
\begin{equation}\label{rrate}
\wh{I}(x^{1},...,x^{k})=\sum_{i=1}^{k}I(x^{i}),
\end{equation}
where
\begin{equation}\label{rate} I(x^{i})=\int_{0}^{\infty}
\sup_{\theta<\theta_{+}}\{\theta\dot
x^{i}(t)-\lambda[\varphi(\theta)-1] \} \rd t,
\end{equation}
if $x^{i}$ are the absolutely continuous functions. We do not give
a definition of functional on  not  absolutely continuous
functions because  condition \reff{phi+} permits to avoid them.

The event $\mathcal{F}$ where minimum of rate function we have to
find is defined by the variational problem presented below.

First let us extend the mapping \reff{jumpo} onto the whole set
$\wh{\mathcal{X}}$ of not decreasing functions
\begin{equation}\label{G}
G:\:\wh{\mathcal{X}}\rightarrow\wh{\mathcal{X}}.
\end{equation}
The extension is denoted by the same symbol.

Consider $X=(x^1,...,x^k)\in \wh{\mathcal{X}}$ and the functions
$\alpha^{i}(t)\in [0,1]$,\, $i=1,...,k$. Let
\[y^{i}(t)=\alpha^{i}(t)x^{i}(t)-(1-\alpha^{i+1}(t))x^{i+1}(t),
\,i=1,...,k.\] Here $Y=GX=(y^{1},...,y^{k})$ is a solution to the
following optimization problem
\begin{equation}\label{opt}
E=\min_{\ovl{\alpha}}\sum_{i}\int|y^{i}(t)-y^{i+1}(t)|{\rm d}t,
\end{equation}
 $\ovl{\alpha}=(\alpha^{1}(t),...,\alpha^{k}(t))$.

The event is
\begin{equation}\label{FF}
\mathcal{F}=\left\{X=(x^{1},...,x^{k}):\:\,\min\left\{\sup_{t\geq 0
}\big\{{y}^{1}(t)-t\big\},\;\sup_{t\geq 0
}\big\{{y}^{k}(t)-t\big\}\right\}>d,\ Y=GX\right\}.
\end{equation}

It follows from the routing rules that the flows  $(y^{i})$ upon
the servers given the input flows $(x^{i})$  are defined by
solution to \reff{opt}. That is the consequence of the system
routing rules. We omit the proof of this fact.

To prove the theorem one needs to find the minimum
\begin{equation}\label{minmin}
\wh{I}(\mathcal{F})=\min\left\{\wh{I}(X):\:X\in \mathcal{F}\right\},
\end{equation}

As the rate function is an integral functional on
$\wh{\mathcal{X}}$, the solution to \reff{minmin} is reduced to
the solution of variational problem on the set of non-decreasing
functions from $\mathcal{F}$.

In fact it is sufficient to find \reff{minmin}  on the piecewise
linear functions of form \reff{xx}. That is because: 1) condition
\reff{phi+} permits to restrict oneself by  absolutely continuous
functions; and 2)  for the homogeneous processes with independent
increments the rate function on the trajectories that connect two
fixed points during a given period takes its minimum on the line
that connects these points. We omit the detailed explanation of
these facts (see \cite{LS,DP2}).

Below we use  the following
 not   formal terminology. The
trajectory $X=(x^1,...,x^k)$ is called {\it an input flow}  and the
image $Y=GX=(y^1,...,y^k)$ {\it a load flow }. A small
neighborhood of trajectory $X$ contents "the real" (jump-wise)
trajectories of input processes $\zeta^i_n$, and a  small
neighborhood of trajectory $Y$ contents "the real" trajectories of
load processes $\hat{w}^i_n$.

  A trajectory $X=(x^1,...,x^k)\in {\cal X}$  that represents  $k$
  input flows will be called a {\it input configuration } or
  simpler a {\it  configuration }
We consider only trajectories of form \reff{xx}, each is
characterized by a pair $(a^{i},T^{i})$ where $a^{i}=\dot x^{i}(t)$
on the interval $[0,T^{i}]$. We denote by $X=\{\ovl{a},T\}$ a
configuration where all $T^{i}$  are equal, here
$\ovl{a}=(a^{1},...,a^{k})$.

For functions of form \reff{xx} a solution to \reff{opt} is
equivalent to solution to  following optimization problem
\begin{equation*}
\hspace{-5cm}\mathbf{O1}\hspace{5cm}
D=\min_{\ovl{\alpha}}\sum_{i=1}^{k}\Big|b^{i}-b^{i+1}\Big|,
\end{equation*}
where
 $$b^{i}=\alpha^{i}a^{i}+(1-\alpha^{i+1})a^{i+1},$$
\begin{equation*}
\alpha^{i}\in [0,1].
\end{equation*}
$\ovl{\alpha}=(\alpha^{1},...,\alpha^{k})$ ($i+1$ mod $k$).

Obviously the solution  exists.

 In our problem for given input flow configuration
 $X=\{\ovl{a},T\}$ vector $\overline b$  represents the load
configuration $Y=\{\ovl{b},T\}$.

In addition to the  above correspondence of all $X$ with $Y$ we
consider also the correspondence of subset of input  flows with a
subset of load flows to the assigned  servers.

We are interested only in {\it connected} subsets
$F^{'}=\{f_{r+1},...,f_{r+l}\}\subseteq F$ of input flows. Let
$X^{'}=\{x^{r+1},...,x^{r+l}\}$ be a configuration of such  flows.
It is supposed that $x^{i}\in X^{'} $ are of  form \reff{xx}, all
$T$ are equal, i.e. $X^{'}=\{\ovl{a}_{r+1,r+l},T\}$ where
$\ovl{a}_{r+1,r+l}=(a^{r+1},...,a^{r+l})\in \mathbb{R}^{l}_{+}$
and $a^{i}=\dot x^{i}(t), t\leq T$. The corresponding
configuration of load flows to the assigned servers
$S(X^{'})=\{s_{r},s_{r+1},...,s_{r+l}\}$ is
$Y^{'}=\{\ovl{b}_{r,r+l},T\}=\{(b^{r},b^{r+1},...,b^{r+l}),T\} $.
Here vector $\ovl{b}_{r,r+l}\in\mathbb{R}^{l+1}$ is defined by a
solution to   optimization problem
\begin{equation*}\hspace{-5cm}\mathbf{O2}\hspace{5cm}
D_{r,r+l}=\min_{\ovl{\alpha}_{r,r+l}}
\sum_{i=r+1}^{r+l}\Big|b^{i-1}-b^{i}\Big|,
\end{equation*}
where
\begin{equation*}
b^{i}=\begin{cases} (1-\alpha^{i+1})a^{i+1}&\mbox{ if }i=r,\\
(1-\alpha^{i+1})a^{i+1}+\alpha^{i}a^{i}&\mbox{ if }r< i\leq
r+l-1,\\
\alpha^{i}a^{i}&\mbox{ if }i= r+l,
\end{cases}
\end{equation*}
\begin{equation*}
\alpha^{i}\in [0,1],
\end{equation*}
$\ovl{\alpha}_{r,r+l}=(\alpha^{r},...,\alpha^{r+l})$.

We look now for solutions to $ \mathbf{O2}$ in case $D_{r,r+l}=0$.
Below the notation  $D_{X,F^{'}}$ is used to indicate that $D$ is
calculated with respect to a set of input flows $F^{'}$ for a
given configuration of all flows $X$. A connected set $F^{'}$ is
called \textit{balanced} if $D_{X,F^{'}}=0$.

  As we are interested in large deviations keeping  in mind $f_1$  only
 subsets $F^{'}\subseteq F$ with $f_{1}\in F^{'}$ are considered.

A balanced $F^{'}\subseteq F$ is said to be \textit{maximal} if
$D_{X,F^{''}}>0$ for any $F^{''}\supset F^{'}$.  A configuration
$X=\{\ovl{a},T\}$ may posses several maximal balanced subsets.

Our aim is to find for a configuration $X=\{\ovl{a},T\}\in {
\mathcal{F}}$ for which the rate function is minimal.

 The rate function for $X=\{\ovl{a},T\}$ is equal to a sum
\begin{equation}\label{inputrate}
I(X)=T\sum_{i=1}^{k}
\big(\theta_{i}a^{i}-\lambda[\varphi(\theta_{i})-1]\big),
\end{equation}
where  $\theta_{i}$ are defined by
$a^{i}=\lambda\varphi^{\prime}(\theta_{i})$.

 Let
$X=\{\ovl{a},T\}$ be a configuration,  $ F^{'}$  a subset of input
flows and $\mathcal{Z}$ a set of configurations such that
\begin{equation*}\label{conf}
\mathcal{Z}( F^{'})=\{Z=\{\ovl{c},T\}:\: c^{i}=a^{i},\mbox{ if
}f_{i}\in F^{'}\}.
\end{equation*}
Rate function for this set is
\begin{equation}\label{balance}
I(\mathcal{Z}( F^{'}))=I(X_0)=T\sum_{i:\:f_{i}\in F^{'}}
\big(\theta_{i}a^{i}-\lambda[\varphi(\theta_{i})-1]\big).
\end{equation}
where $X_0=\{\ovl{a},T\}$ is a configuration with $a_0^i=a_i$ if
$f_i\in F'$ and $a_0^j=a^{(0)}$ if $f_j\not\in F'$


For balanced with respect to $X=\{\ovl{a},T\}$  subset $F^{'}$ we
consider
\begin{equation*}
h_{X,F^{'}}=\frac{\sum_{i:\:f_{i}\in F^{'}}a^{i}}{|F^{'}|}.
\end{equation*}
Here $|F^{'}|$ is the number of flows in $F^{'}$.

 Let $h>1$ and let
$\mathcal{Z}(h,F^{'})=\{Z=(\ovl{a},T),\,h_{Z,F^{'}}=h\}$ be a set of
trajectories with fixed  mean value $h$ on $F'$.
\medskip

\begin{lemma}\label{la1}
 For a set of configurations $\mathcal{Z}(h,F^{'})$ a rate function  $I(\mathcal{Z}(h,F^{'}))$
  on interval $[0,T]$ is equal to
\begin{equation}\label{a1.1}
I(\mathcal{Z}(h,F^{'})=\left(h\wt{\theta}-\lambda[\varphi(\wt{\theta})-1]\right)
T|F^{'}|,
\end{equation}
where $\wt{\theta}$ is determined by equality
\begin{equation}\label{a1h}
h=\lambda\varphi^{\prime}(\wt{\theta}).
\end{equation}

 If $F^{'}$ is a maximal balanced set then
\begin{equation}\label{hh}
\lim_{\varepsilon \to 0}\lim_{n\to
\infty}\frac{-1}{n}\ln\Pr\left(\mathcal{A}_{\varepsilon,n}(h,F^{'})\right)=
I(\mathcal{Z}(h,F^{'})).
\end{equation}
Here we denote
 $$
\mathcal{A}_{\varepsilon,n}(h,F^{'})=\bigcup_{Z\in\mathcal{Z}(h,F^{'})}
\left[ \bigcap_{i:\:f_{i}\in
F^{'}}\mathcal{A}_{\varepsilon,n}^{i}(z_{T,c^{i}}^{i})
\cap\bigcap_{j:\:f_{j}\in F\setminus
F^{'}}{\mathcal{A}}_{\varepsilon,n}^{j}(z_{T,a^{(0)}}^{j})\right],
$$
 $Z=\{\ovl{c},T\}=\{z^{1},...,z^{k}\}$.
\end{lemma}

This Lemma shows that equal overheating in all  flows from $F'$ is
"more" probable than the not equal one. We write "more" because
that is an asymptotic result. In fact the considered
probabilities decay exponentially and the exponent determined by
rate function $I$  is minimal in case of equal flows.

\medskip

\proofr Rewrite the expression for \reff{inputrate}
\begin{equation*}
I({\cal Z}(h,F'))=\sum_{i:\:f_{i}\in
F^{'}}\big((h+g_{i})\theta_{i}-\lambda[\varphi(\theta_{i})-1]\big)T,
\end{equation*}
where $g_{i}=c^{i}-h$. The solution to  system
 $\frac{\partial I({\cal Z})}{\partial g_{i}}=0$ is:
$\theta_{i}=\theta_{j}$ for all $i,j:\:f_{i},f_{j}\in F^{'}$. The
solution is unique thanks to monotonicity in $\theta$ of
$\frac{\partial (c\theta -\lambda(\varphi(\theta)-1))}{\partial
\theta}$, $c=\lambda \varphi'(\theta)$. Equality \reff{hh}  followes
from \reff{balance}.

$\blacktriangle$

\medskip
It follows from  $X\in {\cal F}$ that
$$bT-T\geq d,
$$
where $b=min\{ b^k,b^{1}\}$ and $\ovl{b}$ is a solution to ${\bf
O1}$.

\medskip

Introduce now a set of configurations  $\mathcal{C}(h,F^{'})$ where
$F^{'}$ is  a single maximal balanced set for $Z=\{\ovl{c},T\}\in
\mathcal{C}(h,F^{'})$, $c^{i}=h$, $f_{i}\in F^{'}$, $h>1$.
 Let
\begin{equation}\label{CC}
\mathcal{C}=\bigcup_{F^{'}\subseteq
F}\bigcup_{h>1}\bigcup_{T>0}\mathcal{C}(h,F^{'}).
\end{equation}
We want to find the rate function for $\mathcal{C}$.
 To find
$I(\mathcal{C})$ as $\lambda$ is fixed one has to minimize ${\cal
C}$ in $T,\,h$ and  $F^{'}$ (see \reff{a1.1}, \reff{a1h}). Because
of circular structure of the network instead of minimization in
$F^{'}$ one can minimize in  number  $l=|F^{'} |$ of input flows
$F^{'}$. If $S'$ is the set of servers assigned to $F'$ then
$|S'|=l+1$ as $l<k$ and  $|S'|=k$ as $l=k$.

Below the calculations are based on the following argument:
the overheat of flows $F'$ that bring the overload of
servers $S'$ can be considered in case $l<k$ as a overload
 of a one-channel system with a server of speed $l+1$ and input
flow of rate $\lambda l$. In case $l=k$ a one-channel system
has a server speed $k$  and the    input flow rate  $\lambda k$.
We call such a one-channel system an {\it auxiliary} system.

Consider  $\mathcal{C}(F^{'})=
\bigcup_{h>1}\bigcup_{T>0}\mathcal{C}(h,F^{'})$ and look for
$I(\mathcal{C}(F^{'}))$ as $|F^{'}|$ is fixed.

Suppose first that $l=|F^{'}|< k-1$. Obviously $F^{'}$ is assigned
to   $|S^{'}|=l+1$ servers. We can consider only such load flows
that form configuration $Y=\{\ovl{b},T\}$, with
$b^{i}=b=\frac{l}{l+1}h$, $s_{i}\in S^{'}$. This configuration and
its small neighborhood belongs to the manifold ${\cal F}$ if
\begin{equation*}
\frac{l}{l+1}hT-T\geq d
\end{equation*}
Now by  \reff{a1.1} and  \reff{a1h} we get
\begin{equation}\label{TT}
I(\mathcal{C}(h,F^{'})=\inf_{T}\left(\wt{\theta}\frac{(d+T)(l+1)}{Tl}
-\lambda[\varphi(\wt{\theta})-1]\right)lT
\end{equation}
where $\wt{\theta}$ is defined by
$\displaystyle{\frac{l(d+T)}{T(l+1)}
=\lambda\varphi^{'}(\wt{\theta})}$. It is easy to see that
\reff{TT} has its infinum as
\begin{equation}\label{2;1}
(l+1)\wt{\theta}=l\lambda[\varphi(\wt{\theta})-1].
\end{equation}
Denote by $\theta(\lambda,l)$ a positive solution to the last
equation. Then we get that
\begin{equation}\label{2;2}
J(\lambda,l)=\inf_{T}I(\mathcal{C}(F^{'})))=(l+1)\theta(\lambda,l)d
\end{equation}
and the optimal $b$ and $T$ are
\begin{equation}\label{btl}
  b=b(l) = \lambda \frac{l}{l+1}
\varphi^{\prime}(\theta(\lambda,l)), \ \ \  T=T(l)
=\frac{d}{\lambda\frac{l}{l+1}\varphi^{\prime}(\theta(\lambda,l))-1}.
\end{equation}
In case where  $l=k$ the sum of servers speed is $k$ and the sum
of input flow rates is $\lambda k$. The configuration
$Y=\{\ovl{b},T\}$ of flows to the servers is such that
$b^{i}=b=h$ for $s_{i}\in S$   and    we have   the  inequality
$hT-T\geq d$. Therefore  \reff{TT} becomes
\begin{equation*}\label{1;11}
I(\mathcal{C}(h,F^{'})) =Tk\left(\wt{\theta}\frac{T+d}{T}-
\lambda[\varphi(\wt{\theta})-1]\right).
\end{equation*}
Optimization in $T$ and $h$ gives
\begin{equation}\label{2;21}
J(\lambda,k)=\inf_{T,h}I(\mathcal{C}(h,F^{'}))=k\theta^{\ast}(\lambda)d,
\end{equation}
where $\theta^{\ast}(\lambda)$ is a positive root to
\begin{equation}\label{2;11}
\wt{\theta}=\lambda[\varphi(\wt{\theta})-1].
\end{equation}
The optimal $b$ and $T$ are
\begin{equation}\label{btk}
 b=b^{*} = \lambda \varphi^{\prime}(\theta^{\ast}(\lambda)), \ \ \
T=T^{*}
=\frac{d}{\lambda\varphi^{\prime}(\theta^{\ast}(\lambda))-1}
\end{equation}

Remark that $J(\lambda,k)<J(\lambda,k-1)$ for any  $\lambda$.
Really, if $|L|=l=k-1$ then all servers are overloaded,
 $J(\lambda,k-1)$ $=$ $k\theta(\lambda,k)d$, where
$\theta(\lambda,k-1)$ is a solution to
${\theta=\lambda\frac{k-1}{k}[\varphi(\theta)-1]}$. At the same time
$J(\lambda,k)=k\theta^{*}(\lambda)d$, where $\theta^{*}(\lambda)$
 is a solution to \reff{2;11}. Therefore
$\theta^{*}(\lambda)<\theta(\lambda,k-1)$, see Fig.2.

To get rate function  $I(\mathcal{C})$ as $k$ and  $\lambda$ are
fixed we have to find  $l,\;l\leq k,$ that brings
$$
\min\Big[ k\theta^*(\lambda),\;\min_{1\leq
l<k}(l+1)\theta(\lambda ,l)\Big].$$

\bigskip

\begin{picture}(460,140)
\put(183,115){$\lambda_1$} \qbezier(110,20)(180,30)(180,110)
\put(210,40){$\lambda_2>\lambda_1$} \put(146,110){$\lambda_2$}
\qbezier(110,20)(160,30)(160,110) \put(110,20){\line(1,0){90}}
\put(200,11){$\theta$} \put(110,20){\line(0,1){110}}
\put(75,130){$\lambda(\varphi(\theta)-1)$}
\put(110,20){\line(1,1){90}} \put(140,1){{\bf Fig. 2.}}
\end{picture}

\medskip


\begin{lemma}\label{le2}
\begin{enumerate}

\item [1)]   For any $k$, $l,\;1\leq l\leq k-1$  there exist such
$\lambda^*_{k,l}$, $0<\lambda^*_{k,l}<\wh{\lambda}$  that
$J(\lambda,l)< J(\lambda,k)$ as $\lambda<\lambda^*_{k,l}$ and
$J(\lambda,l)> J(\lambda,k)$ as $\lambda>\lambda^*_{k,l}$.
\item [2)] For any $l_1$, $l_2,\;1\leq l_1<l_2$  there exist such
$\lambda_{l_2,l_1}$,
$0<\lambda_{l_2,l_1}<\displaystyle{\frac{l_2+1}{l_2}\wh{\lambda}}$
that $J(\lambda,l_1)< J(\lambda,l_2)$ as $\lambda<\lambda_{l_2,l_1}$
and $J(\lambda,l_1)> J(\lambda,l_2)$ as $\lambda>\lambda_{l_2,l_1}$.
\end{enumerate}
\end{lemma}
\medskip

{\proofr} 1).   For the start we find such
$\lambda=\lambda^*_{k,l}$,
 $l\leq k-1$ that $J(\lambda,l)=J(\lambda,k)$. Denote
 $(l+1)\theta=\vartheta_l,\;l<k$,\; $k\theta=\vartheta^*_k$.
 It follows from
  \reff{2;11} and \reff{2;1} that
\begin{equation}\label{vartheta}
\vartheta_k^*=\lambda k\Big[\varphi\Big(\frac{\vartheta_k^*}{k}\Big)-1\Big],\;\;\;
\vartheta_l=l\lambda\Big[\varphi\Big(\frac{\vartheta_l}{l+1}\Big)-1\Big].
\end{equation}
By \reff{2;2} and \reff{2;21} the needed equality is achieved as
$\vartheta_l=\vartheta_k^*$. Let us show  that equation
\begin{equation}\label{thetakl}
l\Big[\varphi\Big( \frac{\vartheta}{l+1}\Big) -1\Big]= k\Big[
\varphi\Big( \frac{\vartheta}{k}\Big) -1\Big]
\end{equation}
has a unique solution  $\vartheta^*_{k,l}$ and
\begin{equation*}
 0<\lambda^*_{k,l}=\frac{\vartheta^*_{k,l}}{k(\varphi(\vartheta^*_{k,l}/k)-1)}<\hat{\lambda}.
 \end{equation*}
 Really, $\varphi(0)=1$, therefore
$$\lim_{\vartheta\to 0}\frac{l\Big[ \varphi \Big( \frac{\vartheta}{l+1}\Big)-1\Big]}
{k\Big[ \varphi\Big(\frac{\vartheta}{k}\Big)-1\Big]}
=\frac{l}{l+1}<1.$$
Further, if $\theta_+<\infty$ then
 $$\lim_{\vartheta/(l+1) \to \theta_+}
  \varphi\Big(\frac{\vartheta}{l+1}\Big)=\infty,
\;\;\;\lim_{\vartheta /(l+1)\to\theta_+}\varphi\Big(\frac{\vartheta}{k}\Big)<\infty.$$
 And if $\theta_+=\infty$ then
$$\lim_{\vartheta \to\infty}\frac{\varphi\Big(\frac{\vartheta}{l+1}\Big)}
{\varphi\Big(\frac{\vartheta}{k}\Big)}=\infty.$$
 That means that there exists
$\vartheta=\vartheta^*_{k,l}$ such that  (\ref{thetakl}) takes
place. The uniqueness follows from because $\varphi$ and
 all its derivatives are convex..

The function $\theta(\lambda)$, $0<\theta<\theta_+$ presented  by
\reff{2;11} is defined for $0<\lambda<\hat{\lambda}$ and
$\theta(\lambda)$ monotonically decreases in  $\lambda$, (see Fig.
2).  Therefore there exists such $\lambda^*_{k,l}$,
$\lambda^*_{k,l}<\hat{\lambda}$, that corresponds to
$\vartheta=\vartheta^*_{k,l}$ and for which the conditions of Lemma
are fulfilled.

\medskip

2). The proof of this item follows once again from existence and
uniqueness of solution $\vartheta_{l_2,l_1}$ to
\begin{equation*}
l_1\Big[\varphi\Big( \frac{\vartheta}{l_1+1}\Big) -1\Big]= l_2\Big[
\varphi\Big( \frac{\vartheta}{l_2+1}\Big) -1\Big],
\end{equation*}
therefore it repeats the proof of item 1). But here we have to
notice that $\theta(\lambda,l)$, $0<\theta<\theta_+$, presented by
\reff{2;1} is defined for
$\displaystyle{\lambda<\hat{\lambda}\frac{l+1}{l}}$ (in \reff{2;11}
we had $\lambda<\hat{\lambda}$). Thus
$0<\lambda_{l_2,l_1}<\hat{\lambda}\frac{l_2+1}{l_2}$ and in general
case it may be that $\lambda_{l_2,l_1}>\hat{\lambda}$ (see Fig. 3).\
\
\ $\blacktriangle$
\medskip

The proof of Theorems 1 and 2 follows from the Lemmas.

 Let us set
$$\lambda_k=\min_{l<k}[\lambda^*_{k,l},\lambda_{l,1}],$$
$$\lambda^k=\max_{l}\lambda^*_{k,l}.$$
By Lemma 2 $J(\lambda,1)<\displaystyle{\min_{1<l<k}J(\lambda,l)}$ as
$\lambda<\min_{l}\lambda_{l,1}$ and therefore for fixed $k$ we have
$J(\lambda,1)<\min_{1<l<k}
\{J(\lambda,l),J(\lambda,k)\}$ as $\lambda<\lambda_k$. The way
we get  \reff{2;2} indicates that as $\lambda<\lambda_k$ the
statements  1) and 2) of the Theorem take place  and the values
\reff{btl} correspond to the values of idem 3) of the Theorem.

Further, $J(\lambda,k)<\displaystyle{\min_{1\leq l<k}J(\lambda,l)}$
as $\lambda>\lambda^k$. The way we get  \reff{2;21}
 indicates that as $\lambda>\lambda^k$ the statements 4) and
5) of the Theorem  2 take place and the values \reff{btk} coincide with
the values of item  6) of the Theorem 2.
 \ \ $\blacktriangle$

\bigskip

\section{Different distributions of message length. Examples}

It is clear that if $\lambda_k=\lambda^k$ for some $k$ then as
 $f_1$ is overheated then, depending on $\lambda$, most probably
either
 only $f_1$ or all flows are overheated. The questions are: when
$\lambda_k=\lambda^k$, what happens if $\lambda_k<\lambda^k$ ?

To answer these questions we look at the location of curves
$\vartheta_l$ and $\vartheta_k^*$ (see \reff{vartheta}) on
$(\lambda,\vartheta)$ plane.

\begin{utv}
If $\lambda_{l,1}>\hat{\lambda}$ for any  $1<l<k$, then
$\lambda_k=\lambda^k$,
\ $\lambda^k\to \hat{\lambda}$ as  $k\to \infty$.
\end{utv}

{\proofr} By Lemma 2  each pair of curves
$\vartheta_{1}(\lambda)$,\;$\vartheta_{l}(\lambda)$ defined by,
 \reff{vartheta} has a unique  point of intersection
 $(\lambda_{l,1},\vartheta_{l,1})$ and
$\vartheta_1<\vartheta_l,\;l>1$, as $\lambda<\lambda_{l,1}$.
Therefore from $\lambda_{l,1}>\hat{\lambda}$ follows  that
$J(\lambda,1)<J(\lambda,l)$ for all $\lambda<\hat{\lambda}$. Thus
 $\displaystyle{\min_{1\leq l<k}J(\lambda,l)=\min[J(\lambda,1),J(\lambda,k)]}$,
i.e. $\lambda_k=\lambda^k$. The curves $\vartheta_l(\lambda)$ do not
depend on $k$, the curves $\vartheta^*_k(\lambda)$ increase with
$k$, therefore $\lambda^*_{k,1}=\lambda^k \to \hat{\lambda}$ \ (see
Fig.3).\
\ $\blacktriangle$

\begin{picture}(450,160)
\put(1,30){\line(1,0){125}}
\put(1,30){\line(0,1){110}}
\qbezier[80](60,30)(60,55)(60,130)
\put(120,30){\line(-2,1){118}}
\put(90,30){\line(-1,1){88}}
\put(60,30){\line(-1,2){50}}
\put(88,45){\footnotesize{$\vartheta_1$}}
\put(73,35){\footnotesize{$\vartheta_2$}}
\put(17,118){\footnotesize{$\vartheta_k^*$}}
\qbezier(320,30)(240,40)(210,130)
\qbezier(290,30)(240,70)(230,130)
\qbezier[200](260,30)(240,50)(235,130)
\put(301,34){\footnotesize{$\vartheta_1$}}
\put(271,34){\footnotesize{$\vartheta_2$}}
\put(200,30){\line(1,0){125}} \put(200,30){\line(0,1){110}}
\put(200,140){\footnotesize{$\vartheta$}}
\put(1,140){\footnotesize{$\vartheta$}}
\put(238,122){\footnotesize{$\vartheta_k^*$}} \put(1,18){$0$}
\put(60,18){\footnotesize{$\hat\lambda$}}
\put(117,18){\footnotesize{$2\hat\lambda$}} \put(200,18){$0$}
\put(240,10){\footnotesize{\bf b)}} \put(80,10){\footnotesize{\bf
a)}} \put(260,18){\footnotesize{$\hat\lambda$}}
\put(317,18){\footnotesize{$2\hat\lambda$}}
\qbezier[60](260,30)(260,75)(260,130) \put(143,5){{\bf Fig. 3}}
\end{picture}

\medskip

\begin{utv}
If $\vartheta_1(\lambda_0)>\vartheta_l(\lambda_0)$ for some
$l,\;l>1$ and $\lambda_0\leq \hat{\lambda}$ then
$\lambda_k<\lambda^k$ for sufficiently large $k$  and $\lambda^k\to
\hat{\lambda}$ as $k\to \infty$.
\end{utv}

\medskip

{\proofr} It follows from the condition of Lemma that
$\lambda_{l,1}\leq
\hat{\lambda}$. Remember that  $\vartheta_k^*(\lambda)$
 increases in
$k$. If $k$ is sufficiently large then the pairs
  $\vartheta_k^*$, $\vartheta_1$
and $\vartheta_k^*$, $\vartheta_l$ intersect at  points
$\lambda^*_{k,1}$ and  $\lambda^*_{k,l}$ where
$\lambda^*_{k,l}>\lambda_{l,1}$,\ $\lambda^*_{k,1}>\lambda_{l,1}$,
and $\lambda^*_{k,l}>\lambda^*_{k,1}$. Therefore
$\vartheta_k^*(\lambda)>\vartheta_1(\lambda)>\vartheta_l(\lambda)$
as $\lambda_{l,1}<\lambda<\lambda^*_{k,1}$ (see Fig. 4), i.e.
$\min_{1<l<k}J(\lambda,l)<\min[J(\lambda,1),J(\lambda,k]$ for such
$\lambda$ and that means $\lambda_k<\lambda^k$. As in Preposition 1
$\lambda^k\to\hat{\lambda}$ as $k\to \infty$.
 $\blacktriangle$

\medskip

\begin{picture}(320,230)
\put(100,60){\line(1,0){235}}
\put(100,60){\line(0,1){130}}
\put(100,48){$0$}
\qbezier[80](200,60)(200,80)(200,175)
\qbezier[50](123,60)(123,65)(123,145)
\put(123,60){\circle*{1}}
\qbezier[30](153,60)(153,65)(153,120)
\put(153,60){\circle*{2}}
\qbezier[30](175,60)(175,65)(175,124)
\put(175,60){\circle*{2}}
\qbezier[15](179,60)(179,75)(179,91)
\put(179,60){\circle*{1}}
\put(118,48){\footnotesize{$\lambda^*_{k_1,1}$}}
\put(147,48){\footnotesize{$\lambda_{l,1}$}}
\put(167,48){\footnotesize{$\lambda^*_{k,1}$}}
\put(200,48){\footnotesize{$\hat\lambda$}}
\put(320,62){\footnotesize{$\lambda$}}
\put(295,48){\footnotesize{$2\hat\lambda$}}
\qbezier(300,60)(145,90)(112,180)
\qbezier(250,60) (145,90)(130,180)
\qbezier(200,60)(175,75)(165,180)
\qbezier[170](200,60)(135,65)(120,180)
\put(207,67){\footnotesize{$\vartheta_l$}}
\put(275,65){\footnotesize{$\vartheta_1$}}
\put(101,185){\footnotesize{$\vartheta$}}
\put(150,175){\footnotesize{$\vartheta_{k}^*$}}
\put(130,97){\footnotesize{$\vartheta_{k_1}^*$}}
\put(17,30){\footnotesize{Here $k>k_1$,\;\;\;\;
$\lambda_{k_1}=\lambda^{k_1}=\lambda^*_{k_1,1}$,}}
\put(172,30){\footnotesize{ $\lambda_{k}<\lambda^{k}$,}}
\put(225,30){\footnotesize{ $\vartheta_l<\vartheta_1$ \; as\;
$\lambda_{l,1}<\lambda<\lambda^*_{k,1}$}} \put(180,1){\rm {\bf
Fig. 4}}
\end{picture}

\medskip

Below we present several examples to show  the realization of
described scenarios.

\vspace{1cm}

{\bf  1}.  Exponential distribution of message length.
\medskip

  The density of message length distribution is
 $ce^{-cx}$, $c>0$ and
$\varphi(\theta)=\displaystyle{\frac{c}{c-\theta}}$, the condition
of stability is $\lambda<\hat{\lambda}=c$. The equations
\reff{vartheta} are linear having the form
$$\vartheta_l=c+l(c-\lambda),\;\;\vartheta_k^*=k(c-\lambda).$$
It is clear that $\vartheta_1(\lambda)<\vartheta_l(\lambda)$ as
$\lambda<c=\hat{\lambda}$. On  $(\lambda, \vartheta)$ plain all
$\vartheta_l(\lambda)$ intersect at a point $(c,c)$  and
$\lambda_k=\lambda^k=c\frac{k-2}{k-1}$ (See Fig 3\ а\ ).

\medskip

{\bf 2}. The density of message length distribution is
\begin{equation*}
\frac{1}{2}\Big((c+g )e^{-(c+g)x}+(c-g)e^{-(c-g)x}\Big) ,\;\;g<c.
\end{equation*}
Here we have
$\displaystyle{\varphi(\theta)-1
= \frac{\theta (c-\theta)}{(c-\theta)^2-g^2}}.$
The stability condition is
$\lambda<\hat{\lambda}=\displaystyle{\frac{c^2-g^2}{c}}.$
For $\lambda\leq
\hat{\lambda}$ we can find
$\vartheta_l({\lambda})=(l+1)\theta({\lambda})$  using \reff{2;1}
$$(c-\theta)^2-\frac{l}{l+1}{\lambda}(c-\theta) -g^2=0;$$
the needed expression is
$$\vartheta_l({\lambda})=(l+1)\Big[c-\frac{{\lambda}l}{2(l+1)}-
\sqrt{\frac{{\lambda}^2l^2}{4(l+1)^2}+g^2}\ \  \Big]$$
\begin{equation}
\label{varl}
=c+l\Big[c-\frac{\lambda}{2}-\sqrt{\frac{{\lambda}^2}{4}+g^2\frac{(l+1)^2}{l^2}}
\Big].
\end{equation}
  The value in square brackets of \reff{varl} increases in  $l$, thus
$\vartheta_l({\lambda})$  increases in $l$ as $\lambda\leq
\hat{\lambda}$. That means that all $\vartheta_l(\lambda)$ intersect as
$\lambda>\hat{\lambda}$ (see Fig 3 b). By Preposition 1 that
indicates that $\lambda_k=\lambda^k$.

\medskip

This example demonstrates the possibility of the following scenario:
the overheat of $l$ connected flows may bring not only overload of
$l +1$ assigned servers, but also  overload of other servers by not
overheated flows.

For example, let $k=3$, $c=1$, $g=0.5$, here $\hat\lambda=0.75$. The
numerical estimation gives $\lambda^*_{3,1}\sim 0.418$, i.e.
$2\lambda^*_{3,1}>\hat{\lambda}$. That means that  as $\lambda$ is
inside an interval
 $\Big(\frac{\hat{\lambda}}{2},\lambda^*_{3,1}\Big)$ then
 the overheated $f_1$
brings not only the overload of $s_1,\ s_3$ but also the overload of
$s_2$ fed by not overheated $f_2$ and $f_3$. And it is easy to
estimate that  $s_1,\ s_3$ are overloaded with  greater speed than
$s_2$.

\medskip

{\bf  3}.  Constant message length.

  Let $\Pr(\xi=x)=\delta
(1/c)$,  $c>0$. Here $\varphi(\theta)=e^{\theta/c}$. The stability
condition is  $\lambda<c$.

 For $\lambda=\hat{\lambda}=c$ we get by  \reff{vartheta} that
$\displaystyle{\frac{\vartheta_l(c)}{c}}=l(e^{\vartheta_1(c)/c(l+1)}-1)$.
It is easy to check that
$\displaystyle{\frac{\vartheta_1(c)}{c}}>2$, and
$\displaystyle{\lim_{l\to
\infty}\frac{\vartheta_l(c)}{c}=2}$, i.e.
$\vartheta_1(\hat{\lambda})>\vartheta_l(\hat{\lambda})$ as $l$ is
sufficiently large. By Preposition 2
 $\lambda_k<\lambda^k$ as $k$ is sufficiently large.

  The numerical estimates performed for $c=1$ and
 $k\leq 35$ indicate that the behavior of  $\lambda_k$ and
$\lambda^k$  changes as $k$  increases (and so do also the scenario
of $\vartheta_l$ and $\vartheta_k^*$ intersection).

\medskip

 a) Here $\lambda_k= \lambda^k =\lambda^*_{k,1}$ as
 $3\leq k\leq 12$,
$\lambda_{k}$ increases in  $k$.

For example $\lambda_k\sim 0.311$ as $k=3$;\ \  $\lambda_k\sim
0.667$ as $k=5$;\ \ $\lambda_k\sim 0.857$ as $k=10$;\ \
$\lambda_k\sim 0.883$ as $k=12$.

\medskip

 b)  As $k>12$ the value  $\lambda_k$ does not change  and is equal to
  $\lambda_{2,1}\sim 0.888$;  $\lambda^k$ increases in $k$,
$\lambda_k<\lambda^k$.

\medskip

 c) Up to $k=28$\ \ $\displaystyle {\min_lJ(\lambda,l)=J(\lambda,2)}$ as
$\lambda_k<\lambda<\lambda^k$. Therefore most probab two flows are overheated.
 For example $\lambda^k\sim 0.910$ as $k=15$;\ \
$\lambda^k\sim 0.935$ as $k=20$;\ \ $\lambda ^k\sim 0.940$ as
$k=25$.

\medskip

 d) As  $k>29$ in addition to interval
 $(\lambda_k,\lambda_{3,2})$, where
$\min_lJ(\lambda,l)=J(\lambda,2)$,  there appears an interval
 $(\lambda_{3,2},\lambda^{k})$ where
$\displaystyle{\min_lJ(\lambda,l)=J(\lambda,3)}$, and most probable
three flows are overheated. Here $\lambda_{3,2}\sim
0.956$;  $\lambda^k\sim 0.959$  as $k=30$; $\lambda ^k\sim
0.965$ as $k=35$.

\medskip

 These estimates show that nonzero interval $(\lambda_k,\lambda^k)$ is small
  and $\lambda^k$ is close to $\hat{\lambda}$ as  $k>12$.
\smallskip
Presented numerical data and some analytic investigation suggest
 that as  $k\to \infty$\; a "jump"
from $l$ to $k$ (where
 $\displaystyle{l=\arg \min_{1<m\leq k}J(m,\lambda)}$ changes
its value from $l<k$ to $k$) happens at $l=o\Big( k\Big)$.
\medskip

 All numerical data are presented with accuracy 0.0005
 \bigskip

We want to remark that in case of large fluctuations  the collective behavior
of dependent servers may take place for others, not circular networks.

\section{Acknowledgement} N.D.V. thanks V.~Blinovski, K.~Duffi,
S.~Pirogov and Yu.~Suhov for the useful discussions. The work of
E.A.P. was partly supported by Grant RUM1-2693-MO-05 of CRDF.


\bigskip

\begin{thebibliography}{99}

\bibitem{AH}{\it Alanyali M.,  Hajek B.} On large deviations
in load sharing networks // Ann. Appl. Probability. 1998, V. 8, │ 1,
P. 67-97 .

\bibitem{T} {\it Turner S.R.E.}  Large deviations  for
Join the Shortest Queue // Fields Inst. Communications. 2000, V. 28,
P. 95-106.

\bibitem{MT} {\it McDonald D.R. and Turner S.R.E.}
 Resource Pooling in Distributed
 Queueing Networks // Fields Inst. Communications. 2000, V. 28, P.107-131.

\bibitem{FD} {\it Foley  R.D.,  McDonald D.R.} Join the shortest queue:
stability and exact asymptotics // Ann. Appl. Probab. 2001, V. 11, │
3, P. 569-607.


\bibitem{PSV} {\it  Pechersky E.A., Suhov Y.M.,  Vvedenskaya N.D.} Large
deviations in a two-server system with dynamic routing // Tech.
report, Isaac Newton Institute for Math. Sci.2003, preprint
NI03075-IGS.

\bibitem {PSV1}{\it N.D. Vvedenskaya, E.A. Pechersky, Y.M. Suhov}
Large Deviations in Some Queueing Systems// Problems Inform. Transmissions,
 2000, V.. 36, │ 1, P. 42-53.

\bibitem{AP} {\it Aspandijarov S.,  Pechersky E.}, One large deviations problem for
compound Poisson processes  in queuing theory // Markov Processes
and Relat. Fields, 1997, V. 3, │ 3, P. 333-366.

\bibitem{DMPSV} {\it  Duffy K.,  Malone D.,  Pechersky E.,   Suhov Y.,
Vvedenskaya N.} Large deviations provide good approximation to
queueing system with dynamic routing //  2004, Tech. report, Dublin
Insitutute for Advanced Studies.

\bibitem{DPSV} {\it Duffy K.,Pechersky  E.A., Suhov Y.M,  Vvedenskaya N.D.}
Using estimated entropy in a queueing system with dynamic routing
// Markov Process and Related Fields. 2007, V.13, │ 1, T. 57-84.

\bibitem{PV} {\it Puhalskii A.A., Vladimirov A.A.}  A large deviation
principal for join the shortest queue //
Mathematics of Operation Research, V. 32, │ 3, P. 700-710, 2007.

\bibitem{DP}{\it
R.L.Dobrushin, E.A. Pechersky} , Large  deviations for random processes with
independent increments on infinite intervals, Problems Inform. Transmissions,
 V. 34, │ 4, 1998, P. 354-382.

\bibitem{LP} {\it Li Z-H, Pechersky  E.} On large deviations in queuing
systems // Resenhas IME-USP 1999, V. 4, │ 2, P. 163-182.

\bibitem{DP2} {\it Dobrushin R.L., Pechersky E.A.}    Large  deviations  for   tandem queuing
systems, Journal of Applied Mathematics and Stochastic Analysis //
7, 3, 1994, 301-330.

\bibitem{LS} {\it Lynch J., Sethuraman J.} Lagre deviations for processes
with independent increments // Ann. Prob. 1987, 15, 2, 610-627.

\bibitem{B} {\it Borovkov A.A.} Stachastic Processes and Queueing Theory, Springer-Verlag, 1976.

\end{thebibliography}
\end{document}